\documentstyle[leqno]{article}

\newtheorem{th}{Theorem}[section]
\newtheorem{prop}[th]{Proposition}
\newtheorem{cor}[th]{Corollary}
\newcounter{defin}[section]
\renewcommand{\thedefin}{\thesection.\arabic{defin}}
\newcounter{ex}[section]
\renewcommand{\theex}{\thesection.\arabic{ex}}
\newcounter{rem}[section]
\renewcommand{\therem}{\thesection.\arabic{rem}}
\title{Metrical Multi-Time Lagrange Geometry\\
of Physical Fields}
\author{Mircea Neagu}
\date{}
\begin{document}
\maketitle
\begin{abstract}
Section 1 contains some physical and geometrical aspects of the Lagrangian
geometry of physical fields developed by Miron and Anastasiei \cite{7},
which represents the start point in our metrical multi-time Lagrangian approach
of the theory of physical fields. Section 2 exposes a geometrization of a
Kronecker $h$-regular Lagrangian function with partial derivatives of order
one, $L:J^1(T,M)\to R$. This geometrization relies on the notion of metrical
multi-time Lagrange space $ML^n_p=(J^1(T,M),L)$ introduced in \cite{12}. We
emphasize that this geometry gives a model for both the gravitational  and
electromagnetic field theory, in a general setting. Thus, Section 3 presents
the metrical multi-time Lagrange theory of electromagnetism and describes
its Maxwell equations. Section 4 presents the Einstein equations which govern
the metrical multi-time Lagrange theory of gravitational field. The conservation
laws of the gravitational field are also described in terms of  metrical
multi-time Lagrange geometry.
\end{abstract}
{\bf Mathematics Subject Classification (1991):} 53C07, 53C43, 53C99\\
{\bf Key words:} 1-jet fibre bundle, metrical multi-time Lagrange space,
Cartan canonical connection, Maxwell equations, Einstein equations.

\section{Lagrangian theory of physical fields}

\hspace{5mm} A lot of geometrical models in Mechanics, Physics or Biology are
based on the notion of ordinary Lagrangian. In this sense, we recall
that  a Lagrange  space $L^n=(M,L(x,y))$ is defined as a pair which consists
of a real, smooth, $n$-dimensional manifold $M$ coordinated by $(x^i)_{i=
\overline{1,n}}$, and a regular Lagrangian
$L:TM\to R$, not necessarily homogenous with respect to the direction $(y^i)_
{i=\overline{1,n}}$. The differential geometry of Lagrange spaces is now
considerably developed and used in various fields to study natural process
where the dependence on position, velocity  or momentum is involved \cite{7}.
Also, the geometry of Lagrange spaces gives a model for both the gravitational
and electromagnetic field, in a very natural blending of the geometrical
structure of the space with the characteristic properties of the physical
fields.

In the sequel, we try to expose the main geometrical and physical aspects of
{\it the Lagrangian theory of physical fields} \cite{7}. In order to do that,
let us consider
\begin{equation}
g_{ij}(x^k,y^k)={1\over 2}{\partial^2L\over\partial y^i\partial y^j},
\end{equation}
the {\it fundamental metrical d-tensor} of an ordinary Lagrangian $L:TM\to R$.
From physical point of view, this d-tensor has the physical meaning of an
{\it "unified"} gravitational field on $TM$, which consists of one {\it
"external"} $(x)$-gravitational field spanned by points $\{x\}$, and the other
{\it "internal"} $(y)$-gravitational field spanned by directions $\{y\}$. It
should be emphasized that $y$ is endowed with some microscopic character of
the space-time structure. Moreover, since $y$ is a vector field different of
an ordinary vector field, the $y$-dependence has combined with the concept of
{\it anisotropy}.

The field theory developed on a Lagrange space $L^n$ relies on a nonlinear
connection $\Gamma=(N^i_j(x,y))$ attached naturally to the given Lagrangian
$L$. This plays the role of mapping operator of internal $(y)$-field on the
external $(x)$-field, and prescribes the {\it "interaction"} between $(x)$-
and $(y)$- fields. From geometrical point of view, the nonlinear connection
allows the construction of the {\it adapted bases}
$
\displaystyle{\left\{{\delta\over\delta x^i}={\partial\over\partial x^i}-
N_i^j{\partial\over\partial y^j},{\partial\over\partial y^i}\right\}\subset
{\cal X}(TM)}
$
and
$
\{dx^i,\delta y^i=dy^i+N^i_jdx^j\}\subset{\cal X}^*(TM).
$

Concerning the {\it "unified"} field $g_{ij}(x,y)$ of $L^n$, the authors
constructed a Sasakian-like metric on $TM$,
\begin{equation}
G=g_{ij}dx^i\otimes dx^j+g_{ij}\delta y^i\otimes\delta y^j.
\end{equation}

As to the spatial structure, the most important thing is to determine the {\it
Cartan canonical connection} $C\Gamma=(L^i_{jk},C^i_{jk})$ with respect to $g_{ij}$,
which  comes from the metrical  conditions
\begin{equation}
\left\{\begin{array}{l}\medskip
\displaystyle{g_{ij\vert k}={\delta g_{ij}\over\delta x^k}-L^m_{ik}g_{mj}-
L^m_{jk}g_{mi}=0}\\
\displaystyle{g_{ij}\vert_k={\partial g_{ij}\over\partial y^k}-C^m_{ik}g_{mj}-
C^m_{jk}g_{mi}=0},
\end{array}\right.
\end{equation}
where $"_{\vert k}"$ and $"\vert_k"$  are  the local $h$-  and $v$- covariant
derivatives of  $C\Gamma$. The importance  to the Cartan canonical connection
comes from its main role played in the Lagrangian theory of physical fields.

In this context, the Einstein equations of the gravitational potentials
$g_{ij}(x,y)$ of a Lagrange space $L^n,\;n>2$, are postulated as being the
Einstein equations attached to $C\Gamma$ and $G$, namely \cite{7},
\begin{equation}
\left\{\begin{array}{ll}\medskip
\displaystyle{R_{ij}-{1\over 2}Rg_{ij}={\cal K}{\cal T}_{ij}^H},&^\prime P_{ij}
={\cal  K}{\cal T}_{ij}^1,\\
\displaystyle{S_{ij}-{1\over 2}Sg_{ij}={\cal K}{\cal T}_{ij}^V},&^{\prime\prime}
P_{ij}=-{\cal  K}{\cal T}_{ij}^2,
\end{array}\right.
\end{equation}
where $R_{ij}=R^m_{ijm}$, $S_{ij}=S^m_{ijm}$, $^\prime P_{ij}=P^m_{ijm}$,
$^{\prime\prime} P_{ij}=P^m_{imj}$ are the Ricci tensors of $C\Gamma$, $R=
g^{ij}R_{ij}$, $S=g^{ij}S_{ij}$ are the scalar curvatures, ${\cal T}_{ij}^H$,
${\cal T}_{ij}^V$, ${\cal T}_{ij}^1$, ${\cal T}_{ij}^2$ are the components of the
energy-momentum tensor ${\cal T}$ and ${\cal K}$ is the Einstein constant (equal to 0 for
vacuum). Moreover, the energy-momentum tensors ${\cal T}_{ij}^H$ and ${\cal T}
_{ij}^V$ satisfy the following {\it conservation laws}
\begin{equation}
{\cal K}{\cal T}^{H\;m}_{\quad j\vert m}=-{1\over 2}(P^{hm}_{js}R^s_{hm}+2R^s_{mj}
P^m_s),\quad{\cal K}{\cal T}^{V\;m}_{\quad j\vert m}=0,
\end{equation}
where all notations are described in \cite{7}.

The Lagrangian theory of electromagnetism relies on the {\it canonical Liouville
vector field} {\bf C}=$\displaystyle{y^i{\partial\over\partial y^i}}$ and the Cartan canonical
connection $C\Gamma$ of the Lagrange space $L^n$. In this context, the authors introduce
the {\it electromagnetic 2-form} on $TM$,
\begin{equation}
F=F_{ij}\delta y^i\wedge  dx^j+f_{ij}\delta y^i\wedge\delta y^j,
\end{equation}
where
\begin{equation}
\begin{array}{l}\medskip
\displaystyle{F_{ij}={1\over 2}[(g_{im}y^m)_{\vert j}-(g_{jm}y^m)_{\vert i}],}\\
\displaystyle{f_{ij}={1\over 2}[(g_{im}y^m)\vert_j-(g_{jm}y^m)\vert_i].}
\end{array}
\end{equation}

Using geometrical identities,  they deduce that the vertical electromagnetic
components $f_{ij}$ vanish always.

At the same time, using the Bianchi identities attached to the Cartan canonical
connection $C\Gamma$, they conclude that the  horizontal electromagnetic
components $F_{ij}$ are governed by the following {\it equations of Maxwell
type},
\begin{equation}
\left\{\begin{array}{l}\medskip
F_{ij\vert k}+F_{jk\vert i}+F_{ki\vert j}=-\sum_{\{i,j,k\}}C_{imr}R^r_{jk}y^m\\
F_{ij}\vert_k+F_{jk}\vert_i+F_{ki}\vert_j=0.
\end{array}\right.
\end{equation}

Finally, we point out that physical aspects of the Lagrangian electromagnetism
are studied by Ikeda in \cite{5}.

In this paper, we naturally extend the
previous field theory to a general one, constructed on the  jet fibre bundle
of order one $J^1(T,M)\to T\times M$, where $T$ is  a smooth, real,  $p$-dimensional
{\it "multi-time"} manifold  coordinated  by $(t^\alpha)_{\alpha=\overline{1,p}}$
and $M$ is a smooth, real $n$-dimensional {\it "spatial"} manifold coordinated
by $(x^i)_{i=\overline{1,n}}$. The gauge group of $J^1(T,M)$ is
\begin{equation}\label{gg}
\left\{\begin{array}{l}\medskip
\tilde t^\alpha=\tilde t^\alpha(t^\beta)\\\medskip
\tilde x^i=\tilde x^i(x^j)\\
\displaystyle{\tilde x^i_\alpha={\partial\tilde x^i\over\partial x^j}
{\partial t^\beta\over\partial\tilde t^\alpha}x^j_\beta}.
\end{array}\right.
\end{equation}
In other words, it is more general than that used in the papers  \cite{7}, \cite{8}.
We recall that the jet fibre bundle of order
one  is a basic object in the study of classical and quantum field theories.

Our field theory is created, in a natural manner, from a given {\it Kronecker
$h$-regular} Lagrangian  function on $J^1(T,M)$ (i. e.  a smooth function
$L:J^1(T,M)\to R$), and can be called the {\it metrical multi-time Lagrange
theory of physical fields.}

In order to have  a clear exposition of our theory, we point ou that we use the
following three distinct notions:

i) {\it multi-time Lagrangian function} $-$ A smooth function $L:J^1(T,M)\to R$.

ii) {\it multi-time Lagrangian} (Olver's terminology) $-$ A local function ${\cal L}$
on $J^1(T,M)$ which transform by the rule $\tilde{\cal L}={\cal L}\vert\det J\vert$, where
$J$ is the Jacobian matrix of coordinate transformations $t^\alpha=
t^\alpha(\tilde t^\beta)$. If $L$ is a Lagrangian function on
1-jet fibre bundle, then ${\cal L}=L\sqrt{\vert h\vert}$
represent a Lagrangian on $J^1(T,M)$.

iii) {\it multi-time Lagrangian density} (Marsden's terminology) $-$ A smooth map
\linebreak ${\cal D}:J^1(T,M)\to \Lambda^p(T^*T)$. For example, the entity ${\cal D}=
{\cal L}dt^1\wedge dt^2\wedge\ldots\wedge dt^p$, where ${\cal L}$ is a
Lagrangian, represents a Lagrangian density on $J^1(T,M)$.

We emphasize that the construction of a theory of physical fields attached
to a given first order multi-time Lagrangian function was tried unsatisfactory,
again, by Miron, Kirkovits and Anastasiei in \cite{8}. In their opinion, a
such construction must be done on the vector bundle $\oplus_1^pTM\to M$, where
the coordinates of
$\alpha$-th copy of $TM$ are denoted $(x^i,x^i_\alpha)$, and its gauge group
is of the form
\begin{equation}
\left\{\begin{array}{l}\medskip
\tilde x^i=\tilde x^i(x^j)\\
\displaystyle{\tilde x^i_\alpha={\partial\tilde x^i\over\partial x^j}x^j_\alpha}.
\end{array}\right.
\end{equation}
In other words, their gauge group ignores the multi-temporal reparametrizations.
From our point of view this is the first difficulty of their theory. At the same
time, their trial was unsatisfactory because they do not succeded to write the
local expressions of the Bianchi identities of the Cartan canonical connection.
This second difficulty of their theory appeared probably from the very complicated
computations that was involved.

In our paper, we remove these difficulties, using a Kronecker $h$-regular
multi-time Lagrange function on the 1-jet fibre bundle $J^1(T,M)$. These
objects allow  us the writing of the Bianchi identities  of the  Cartan
canonical connection, so necesssary in the description of field equations.
Nevertheless, our theory has also a difficulty, coming from the quite strong
condition of Kronecker $h$-regularity imposed to the multi-time Lagrangian function.
This difficulty will be removed in the paper \cite{9}.

\section{Metrical multi-time Lagrange spaces}

\setcounter{equation}{0}
\hspace{5mm}  Let us consider $T$ (resp. $M$) a {\it "temporal"}  (resp. {\it
"spatial"}) manifold of dimension $p$ (resp. $n$), coordinated by $(t^\alpha)_
{\alpha=\overline{1,p}}$ (resp. $(x^i)_{i=\overline{1,n}}$). Let $E=J^1(T,M)\to
T\times M$ be the jet fibre bundle of order one associated to these manifolds.
The {\it bundle of configuration} $J^1(T,M)$ is coordinated by
$(t^\alpha,x^i,x^i_\alpha)$, where $\alpha=\overline{1,p}$ and $i=\overline{1,n}$.
Note that the terminology used above is justified in \cite{11}.\medskip\\
\addtocounter{rem}{1}
{\bf Remarks \therem} i) Throughout  this paper, the indices $\alpha,\beta,\gamma,
\ldots$ run from $1$ to $p$, and the indices $i,j,k,\ldots$ run from $1$ to $n$.

ii) In the particular case $T=R$ (i. e., the temporal manifold $T$ is the usual
time axis represented by the set of real numbers), the coordinates  $(t^1,x^i,
x^i_1)$ of the 1-jet space  $J^1(R,M)\equiv R\times TM$ are denoted $(t,x^i,y^i)$.
\medskip

We start our study considering  a smooth multi-time Lagrangian function
$L:E\to R$, which is locally expressed by $E\ni(t^\alpha,x^i,x^i_\alpha)\to
L(t^\alpha,x^i,x^i_\alpha)\in R$. The {\it vertical fundamental  metrical
d-tensor} of $L$ is
\begin{equation}
G^{(\alpha)(\beta)}_{(i)(j)}={1\over 2}{\partial^2L\over\partial x^i_\alpha
\partial x^j_\beta}.
\end{equation}

Now, let $h=(h_{\alpha\beta})$ be a fixed semi-Riemannian metric on the temporal
manifold $T$ and $g_{ij}(t^\gamma, x^k, x^k_\gamma)$ be a d-tensor on $E$, symmetric,
of rank $n$, and having a constant signature.\medskip\\
\addtocounter{defin}{1}
{\bf Definition \thedefin} A multi-time Lagrangian function $L:E\to R$ whose vertical
fundamental metrical d-tensor is of the form
\begin{equation}
G^{(\alpha)(\beta)}_{(i)(j)}(t^\gamma,x^k,x^k_\gamma)=h^{\alpha\beta}(t^
\gamma)g_{ij}(t^\gamma,x^k,x^k_\gamma),
\end{equation}
is called a {\it Kronecker $h$-regular multi-time Lagrangian function
with respect to the
temporal semi-Riemannian metric $h=(h_{\alpha\beta})$}.\medskip

In this context, we can introduce the  following\medskip\\
\addtocounter{defin}{1}
{\bf  Definition \thedefin} A pair $ML^n_p=(J^1(T,M),L)$, where $p=\dim T$ and $n=\dim M$,
which consists of the  1-jet  fibre bundle and a
Kronecker $h$-regular multi-time Lagrangian function $L:J^1(T,M)\to R$ is called a {\it
metrical multi-time Lagrange space}.\medskip\\
\addtocounter{rem}{1}
{\bf Remarks \therem} i) In the particular case $(T,h)=(R,\delta)$, a metrical
multi-time Lagrange space is called a {\it  relativistic rheonomic Lagrange
space}  and is denoted $RL^n=(J^1(R,M),L)$.

ii) If  the temporal manifold $T$ is 1-dimensional, then,
via a temporal reparametrization, we have  $J^1(T,M)\equiv J^1(R,M)$.
In other words, a metrical multi-time Lagrangian space having $\dim T=1$ is a
{\it reparametrized relativistic rheonomic Lagrange space}.\medskip\\
\addtocounter{ex}{1}
{\bf Examples \theex} i) Suppose that the spatial manifold $M$ is also endowed with
a semi-Riemannian metric $g=(g_{ij}(x))$. Then, the  multi-time Lagrangian function
\begin{equation}
L_1:J^1(T,M)\to R,\quad
L_1=h^{\alpha\beta}(t)g_{ij}(x)x^i_\alpha x^j_\beta
\end{equation}
is  a Kronecker $h$-regular multi-time Lagrangian function. Consequently,
$ML^n_p=(J^1(T,M),L_1)$ is
a metrical multi-time Lagrange space. We underline that the multi-time Lagrangian
${\cal L}_1=L_1\sqrt{\vert h\vert}$ is exactly the energy multi-time Lagrangian
whose extremals are the harmonic maps between the pseudo-Riemannian manifolds
$(T,h)$ and $(M,g)$ \cite{3}. At the same time, this
multi-time Lagrangian is a basic object in the physical theory of bosonic
strings.

ii) In above notations, taking $U^{(\alpha)}_{(i)}(t,x)$ as a d-tensor field on
$E$ and \linebreak $F:T\times M\to R$ a smooth map, the more general multi-time
Lagrangian function
\begin{equation}
L_2:E\to R,\quad
L_2=h^{\alpha\beta}(t)g_{ij}(x)x^i_\alpha x^j_\beta+U^{(\alpha)}_{(i)}(t,x)x^i_
\alpha+F(t,x)
\end{equation}
is also a Kronecker $h$-regular multi-time Lagrangian. The metrical multi-time
Lagrange space $ML^n_p=(J^1(T,M),L_2)$ is called the {\it autonomous metrical
multi-time Lagrange space  of electrodynamics} because, in the particular case
$(T,h)=(R,\delta)$, we recover the classical Lagrangian space of electrodynamics
\cite{7} which governs the movement law  of a particle placed concomitantly into
a gravitational field and an electromagnetic one. From physical point of view,
the semi-Riemannian metric $h_{\alpha\beta}(t)$ (resp. $g_{ij}(x)$) represents the {\it gravitational
potentials} of the space $T$ (resp. $M$), the d-tensor $U^{(\alpha)}_{(i)}(t,x)$
stands for the {\it electromagnetic potentials} and $F$ is a function which is
called {\it potential function}. The  non-dynamical character of spatial gravitational
potentials $g_{ij}(x)$ motivates us to use the term {\it "autonomous"}.

iii) More general, if  we consider $g_{ij}(t,x)$ a d-tensor field  on $E$,
symmetric,  of rank $n$ and  having constant signature on $E$, we can
define the Kronecker $h$-regular multi-time Lagrangian function
\begin{equation}
L_3:E\to R,\quad
L_3=h^{\alpha\beta}(t)g_{ij}(t,x)x^i_\alpha x^j_\beta+U^{(\alpha)}_{(i)}(t,x)x^i_
\alpha+F(t,x).
\end{equation}
The pair $ML^n_p=(J^1(T,M),L_3)$ is a metrical multi-time Lagrange  space which
is called the
{\it non-autonomous metrical multi-time Lagrange space of electrodynamics}.
Physically, we
remark that the gravitational potentials $g_{ij}(t,x)$ of the spatial manifold
$M$ are dependent of the temporal coordinates $t^\gamma$, emphasizing their
dynamic character.

An important role and, at the same time, an obstruction in the subsequent
development  of the metrical multi-time Lagrangian geometry, is played by
the next \cite{11}
\begin{th} {(characterization of metrical multi-time Lagrange spaces)}\\
If we have $\dim T\geq 2$, then the following statements are equivalent:

i) $L$ is a Kronecker $h$-regular multi-time Lagrangian function on $J^1(T,M)$.

ii) The multi-time Lagrangian function $L$ reduces  to a non-autonomous electrodynamics
multi-time  Lagrangian function, that is,
$$
L=h^{\alpha\beta}(t)g_{ij}(t,x)x^i_\alpha x^j_\beta+U^{(\alpha)}_{(i)}(t,x)x^i_
\alpha+F(t,x).
$$
\end{th}

A direct consequence of the previous characterization theorem is the following
\begin{cor}
The fundamental vertical metrical d-tensor of an arbitrary Kronecker
$h$-regular multi-time Lagrangian function $L$ is of the form
\begin{equation}
G^{(\alpha)(\beta)}_{(i)(j)}={1\over 2}{\partial^2L\over\partial x^i_\alpha
\partial x^j_\beta}=\left\{\begin{array}{ll}\medskip
h^{11}(t)g_{ij}(t,x^k,y^k),&p=1\\
h^{\alpha\beta}(t^\gamma)g_{ij}(t^\gamma,x^k),&p\geq 2,
\end{array}\right.
\end{equation}
where  $p=\dim T$.
\end{cor}
\addtocounter{rem}{1}
{\bf Remarks \therem} i) It is obvious that the preceding theorem is  an obstruction
in the development  of a fertile metrical multi-time Lagrangian geometry. This obstruction will
be removed in a subsequent paper by the introduction of a more general notion,
that of {\it generalized metrical multi-time Lagrange space} [9]. The generalized
metrical multi-time Lagrange geometry and its theory of physical fields are
constructed in \cite{9} using
just a given {$h$-regular fundamental vertical metrical $d$-tensor $G^{(\alpha)(\beta)}_{(i)(j)}$
on the 1-jet space $J^1(T,M)$.

ii) In the case $p=\dim T\geq 2$, the above theorem obliges us to continue the
study of the metrical multi-time Lagrangian space theory, channeling our attention upon the
non-autonomous metrical multi-time Lagrange space of electrodynamics.\medskip

Following the geometrical development from the paper \cite  {12}, the fundamental
vertical metrical d-tensor $G^{(\alpha)(\beta)}_{(i)(j)}$ of the metrical
multi-time Lagrange space $ML^n_p=(J^1(T,M),L)$  induces naturally a {\it
canonical nonlinear connection} $\Gamma=(M^{(i)}_{(\alpha)\beta},N^{(i)}_
{(\alpha)j})$ on $E=J^1(T,M)$.
\begin{th}
The canonical nonlinear connection $\Gamma$  of the metrical multi-time Lagrange
space $ML^n_p=(J^1(T,M),L)$ is defined by the temporal components
\begin{equation}
M^{(i)}_{(\alpha)\beta}=\left\{\begin{array}{ll}\medskip
-H^1_{11}y^i,&p=1\\\medskip
-H^\gamma_{\alpha\beta}x^i_\gamma,&p\geq 2
\end{array}\right.
\end{equation}
and the spatial components
\begin{equation}
N^{(i)}_{(\alpha)j}=\left\{\begin{array}{ll}\medskip
\displaystyle{h_{11}{\partial{\cal G}^i\over\partial y^j}},&p=1\\\medskip
\displaystyle{\Gamma^i_{jk}x^k_\alpha+{g^{ik}\over 2}{\partial g_{jk}\over
\partial t^\alpha}+{g^{ik}\over 4}h_{\alpha\beta}U^{(\beta)}_{(k)j}},&p\geq 2,
\end{array}\right.
\end{equation}
where
\begin{equation}
\begin{array}{l}\medskip
\displaystyle{{\cal G}^i={g^{ik}\over 4}\left({\partial^2L\over\partial x^j
\partial y^k}y^j-{\partial L\over\partial x^k}+{\partial^2L\over\partial t
\partial y^k}+{\partial L\over\partial x^k}H^1_{11}+2h^{11}H^1_{11}g_{kl}y^l
\right)},\\\medskip
\displaystyle{H^\gamma_{\alpha\beta}={h^{\gamma\eta}\over 2}\left(
{\partial h_{\eta\alpha}\over\partial t^\beta}+
{\partial h_{\eta\beta}\over\partial t^\alpha}-
{\partial h_{\alpha\beta}\over\partial t^\eta}\right)},\\\medskip
\displaystyle{\Gamma^i_{jk}={g^{im}\over 2}\left(
{\partial g_{mj}\over\partial x^k}+
{\partial g_{mk}\over\partial x^j}-
{\partial g_{jk}\over\partial x^m}\right)},\\\medskip
\displaystyle{U^{(\beta)}_{(k)j}={\partial U^{(\beta)}_{(k)}\over\partial
x^j}-{\partial U^{(\beta)}_{(j)}\over\partial x^k}}.
\end{array}
\end{equation}
\end{th}
\addtocounter{rem}{1}
{\bf Remarks \therem} i) Considering the particular case $(T,h)=(R,\delta)$, we remark
that the canonical nonlinear connection $\Gamma=(0,N^{(i)}_{(1)j})$  of  the
relativistic rheonomic Lagrange space $RL^n=(J^1(R,M),L)$ reduces to the canonical nonlinear
connection from Miron-Anastasiei theory \cite{7}.

ii) In the case of  an autonomous metrical multi-time Lagrange space of
electrodynamics  (i. e., $g_{ij}(t^\gamma,x^k,x^k_\gamma)=g_{ij}(x^k)$), the generalized Christoffel
symbols $\Gamma^i_{jk}(t^\mu,x^m)$ of the metrical d-tensor $g_{ij}$ reduce
to the classical ones $\gamma^i_{jk}(x^m)$, and the canonical nonlinear
connection becomes $\Gamma=(M^{(i)}_{(\alpha)\beta},N^{(i)}_{(\alpha)j})$,
where
$$
M^{(i)}_{(\alpha)\beta}=\left\{\begin{array}{ll}\medskip
-H^1_{11}y^i,&p=1\\\medskip
-H^\gamma_{\alpha\beta}x^i_\gamma,&p\geq 2
\end{array}\right.\quad\mbox{and}\quad
N^{(i)}_{(\alpha)j}=\left\{\begin{array}{ll}\medskip
\displaystyle{\gamma^i_{jk}y^k+{g^{ik}\over 4}h_{11}U^{(1)}_{(k)j}},&p=1\\\medskip
\displaystyle{\gamma^i_{jk}x^k_\alpha+{g^{ik}\over 4}h_{\alpha\gamma}U^
{(\gamma)}_{(k)j}},&p\geq 2.
\end{array}\right.
$$

The main result of the metrical multi-time Lagrange geometry is the theorem
of existence of the {\it Cartan canonical $h$-normal linear connection} $C\Gamma$ which  allow the
subsequent development of the  {\it metrical  multi-time Lagrangian theory of
physical fields}.
\begin{th}{(of existence and uniqueness of Cartan canonical connection)}
On the metrical multi-time Lagrange space $ML^n_p=(J^1(T,M),L)$ endowed with its canonical
nonlinear connection $\Gamma$, there is a unique $h$-normal $\Gamma$-linear
connection
$$
C\Gamma=(H^\gamma_{\alpha\beta},G^k_{j\gamma},L^i_{jk},C^{i(\gamma)}
_{j(k)})
$$
having the metrical properties\medskip

i) $g_{ij\vert k}=0,\quad g_{ij}\vert^{(\gamma)}_{(k)}=0$,\medskip

ii) $\displaystyle{ G^k_{j\gamma}={g^{ki}\over 2}{\delta g_{ij}\over\delta
t^\gamma},\quad L^k_{ij}=L^k_{ji},\quad
C^{i(\gamma)}_{j(k)}=C^{i(\gamma)}_{k(j)}}$.

Moreover, the coefficients $L^i_{jk}$ and $C^{i(\gamma)}_{j(k)}$ of the
Cartan canonical connection have the expressions \cite{12}
\begin{equation}
\begin{array}{l}
\displaystyle{L^i_{jk}={g^{im}\over 2}\left(
{\delta g_{mj}\over\delta x^k}+
{\delta g_{mk}\over\delta x^j}-
{\delta g_{jk}\over\delta x^m}\right),}\\
\displaystyle{C^{i(\gamma)}_{j(k)}={g^{im}\over 2}\left(
{\partial g_{mj}\over\partial x^k_\gamma}+
{\partial g_{mk}\over\partial x^j_\gamma}-
{\partial g_{jk}\over\partial x^m_\gamma}\right).}
\end{array}
\end{equation}
\end{th}
\addtocounter{rem}{1}
{\bf Remarks \therem} i) In the particular case $(T,h)=(R,\delta)$, the Cartan
canonical $\delta$-normal $\Gamma$-linear connection of the relativistic rheonomic
Lagrange space $RL^n=(J^1(R,M),L)$ reduces to the Cartan canonical connection
used in \cite{7}.

ii) As a rule, the  Cartan canonical connection of a metrical multi-time
Lagrange space $ML^n_p$ verifies also the metrical properties
$$
h_{\alpha\beta/\gamma}=h_{\alpha\beta\vert k}=h_{\alpha\beta}\vert^{(\gamma)}
_{(k)}=0\;\mbox{and}\;g_{ij/\gamma}=0.
$$

iii) In the case $p=\dim T\geq 2$, the coefficients of the Cartan connection
of a metrical multi-time Lagrange space reduce to
$$
\bar G^\gamma_{\alpha\beta}=H^\gamma_{\alpha\beta},\;G^k_{j\gamma}={g^{ki}
\over 2}{\partial g_{ij}\over\partial t^\gamma},\;L^i_{jk}=\Gamma^i_{jk},\;
C^{i(\gamma)}_{j(k)}=0.
$$

\begin{th}
The torsion d-tensor {\bf T} of the Cartan canonical connection of a metrical
multi-time Lagrange space is determined by the local components
\begin{equation}
\begin{tabular}{|c|*{2}{c|}*{2}{c|}*{2}{c|}}
\hline
&\multicolumn{2}{|c}{$h_T$}&
\multicolumn{2}{|c}{$h_M$}&
\multicolumn{2}{|c|}{$v$}\\
\cline{2-7}
&$p=1$&$p\geq 2$&$p=1$&$p\geq 2$&$p=1$&$p\geq 2$\\
\hline
$h_Th_T$&0&0&0&0&0&$R^{(m)}_{(\mu)\alpha\beta}$\\
\hline
$h_Mh_T$&0&0&$T^m_{1j}$&$T^m_{\alpha j}$&$R^{(m)}_{(1)1j}$&
$R^{(m)}_{(\mu)\alpha j}$\\
\hline
$h_Mh_M$&0&0&0&0&$R^{(m)}_{(1)ij}$&$R^{(m)}_{(\mu)ij}$\\
\hline
$vh_T$&0&0&0&0&$P^{(m)\;\;(1)}_{(1)1(j)}$&$P^{(m)\;\;(\beta)}_{(\mu)\alpha(j)}$\\
\hline
$vh_M$&0&0&$P^{m(1)}_{i(j)}$&0&$P^{(m)\;(1)}_{(1)i(j)}$&0\\
\hline
$vv$&0&0&0&0&0&0\\
\hline
\end{tabular}
\end{equation}
where,

i) for $p=\dim T=1$, we have\medskip

\hspace*{25mm}$T^m_{1j}=-G^m_{j1}$, $\;P^{m(1)}_{i(j)}=C^{m(1)}_
{i(j)}$, $\;P^{(m)\;(1)}_{(1)1(j)}=-G^m_{j1}$,\medskip

\hspace*{25mm}$\displaystyle{P^{(m)\;(1)}_{(1)i(j)}={\partial N^{(m)}_{(1)i}\over\partial y^j}
-L^m_{ji}}$, $\quad\displaystyle{{\delta N^{(m)}_{(1)i}\over\delta x^j}-
{\delta N^{(m)}_{(1)j}\over\delta x^i}}$,\medskip

\hspace*{25mm}$\displaystyle{R^{(m)}_{(1)1j}=
-{\partial N^{(m)}_{(1)j}\over\partial t}+H^1_{11}\left[N^{(m)}_{(1)j}-{\partial
N^{(m)}_{(1)j}\over\partial y^k}y^k\right]}$;\medskip

ii) for $p=\dim T\geq 2$, denoting
$$
\begin{array}{l}\medskip
\displaystyle{F^m_{i(\mu)}={g^{mp}\over 2}\left[{\partial g_{pi}\over\partial
t^\mu}+{1\over 2}h_{\mu\beta}U^{(\beta)}_{(p)i}\right]},\\\medskip
\displaystyle{H^\gamma_{\mu\alpha\beta}={\partial H^\gamma_{\mu\alpha}\over
\partial t^\beta}-{\partial H^\gamma_{\mu\beta}\over\partial t^\alpha}+
H^\eta_{\mu\alpha}H^\gamma_{\eta\beta}-H^\eta_{\mu\beta}H^\gamma_{\eta\alpha},}\\
\displaystyle{r^m_{pij}={\partial \Gamma^m_{pi}\over\partial x^j}-
{\partial \Gamma^m_{pj}\over\partial x^i}+\Gamma^k_{pi}\Gamma^m_{kj}-\Gamma^
k_{pj}\Gamma^m_{ki},}
\end{array}
$$
we have
$$
\begin{array}{l}\medskip
T^m_{\alpha j}=-G^m_{j\alpha}$, $\;P^{m\;\;(\beta)}_{(\mu)\alpha
(j)}=-\delta^\beta_\gamma G^m_{j\alpha}$, $\;R^{(m)}_{(\mu)\alpha(j)}=-H^\gamma
_{\mu\alpha\beta}x^m_\gamma,\\\medskip
\displaystyle{R^{(m)}_{(\mu)\alpha j}=
-{\partial N^{(m)}_{(\mu)j}\over\partial t^\alpha}+{g^{mk}\over 2}H^\beta_
{\mu\alpha}\left[{\partial g_{jk}\over\partial t^\beta}+{h_{\beta\gamma}\over
2}U^{(\gamma)}_{(k)j}\right]},\\
\displaystyle{R^{(m)}_{(\mu)ij}=r^m_{ijk}x^k_\mu+\left[
F^m_{i(\mu)\vert j}-F^m_{j(\mu)\vert i}\right]}.
\end{array}
$$
\end{th}
\addtocounter{rem}{1}
{\bf Remark \therem} In the case of autonomous metrical multi-time Lagrange
space of electrodynamics
(i. e., $g_{ij}(t^\gamma,x^k,x^k_\gamma)=g_{ij}(x^k)$), all torsion d-tensors
of the Cartan connection vanish,  except
$$
R^{(m)}_{(\mu)\alpha\beta}=-H^\gamma_{\mu\alpha\beta}x^m_\gamma,\quad
R^{(m)}_{(\mu)\alpha j}=-{h_{\mu\eta}g^{mk}\over 4}\left[H^\eta_{\alpha\gamma}
U^{(\gamma)}_{(k)j}+{\partial U^{(\eta)}_{(k)j}\over\partial t^\alpha}\right],
$$
$$
R^{(m)}_{(\mu)ij}=r^m_{ijk}x^k_\mu+{h_{\mu\eta}g^{mk}\over 4}\left[
U^{(\eta)}_{(k)i\vert j}+U^{(\eta)}_{(k)j\vert i}\right],
$$
where $H^\gamma_{\mu\alpha\beta}$ (resp. $r^m_{ijk}$) are the curvature tensors
of the semi-Riemannian metric $h_{\alpha\beta}$ (resp. $g_{ij}$).
\begin{th}
The  curvature d-tensor {\bf R} of the Cartan
canonical connection is determined by the local components
$$
\begin{tabular}{|c|*{2}{c|}*{2}{c|}*{2}{c|}}
\hline
&\multicolumn{2}{|c}{$h_T$}&
\multicolumn{2}{|c}{$h_M$}&
\multicolumn{2}{|c|}{$v$}\\
\cline{2-7}
&$p=1$&$p\geq 2$&$p=1$&$p\geq 2$&$p=1$&$p\geq 2$\\
\hline
$h_Th_T$&0&$H^\alpha_{\eta\beta\gamma}$&0&$R^l_{i\beta\gamma}$&0&
$R^{(l)(\alpha)}_{(\eta)(i)\beta\gamma}$\\
\hline
$h_Mh_T$&0&0&$R^l_{i1k}$&$R^l_{i\beta k}$&$R^{(l)(1)}_{(1)(i)1k}=R^l_{i1k}$&
$R^{(l)(\alpha)}_{(\eta)(i)\beta k}$\\
\hline
$h_Mh_M$&0&0&$R^l_{ijk}$&$R^l_{ijk}$&$R^{(l)(1)}_{(1)(i)jk}=R^l_{ijk}$&
$R^{(l)(\alpha)}_{(\eta)(i)jk}$\\
\hline
$vh_T$&0&0&$P^{(l)\;\;(1)}_{i1(k)}$&0&$P^{(l)(1)\;(1)}_{(1)(i)1(k)}=
P^{(l)\;\;(1)}_{i1(k)}$&0\\
\hline
$vh_M$&0&0&$P^{l\;(1)}_{ij(k)}$&0&$P^{(l)(1)\;(1)}_{(1)(i)j(k)}=P^{l\;(1)}_
{ij(k)}$&0\\
\hline
$vv$&0&0&$S^{l(1)(1)}_{i(j)(k)}$&0&$S^{(l)(1)(1)(1)}_{(1)(i)(j)(k)}=S^{l(1)(1)}
_{i(j)(k)}$&0\\
\hline
\end{tabular}
$$
where
$R^{(l)(\alpha)}_{(\eta)(i)\beta\gamma}=\delta^\alpha_\eta R^l_{i\beta\gamma}+
\delta^l_i H^\alpha_{\eta\beta\gamma}$,
$R^{(l)(\alpha)}_{(\eta)(i)\beta k}=\delta^\alpha_\eta R^l_{i\beta k}$,
$R^{(l)(\alpha)}_{(\eta)(i)jk}=\delta^\alpha_\eta R^l_{ijk}$ and\medskip

i) for $p=\dim T=1$, we have\medskip

$
\displaystyle{R^l_{i1k}={\delta G^l_{i1}\over\delta x^k}-
{\delta L^l_{ik}\over\delta t}+G^m_{i1}L^l_{mk}-
L^m_{ik}G^l_{m1}+C^{l(1)}_{i(m)}R^{(m)}_{(1)1k},}
$\medskip

$
\displaystyle{R^l_{ijk}={\delta L^l_{ij}\over\delta x^k}-
{\delta L^l_{ik}\over\delta x^j}+L^m_{ij}L^l_{mk}-L^m_{ik}L^l_{mj}+
C^{l(1)}_{i(m)}R^{(m)}_{(1)jk},}
$\medskip

$
\displaystyle{P^{l\;(1)}_{i1(k)}={\partial G^l_{i1}\over\partial
y^k}-C^{l(1)}_{i(k)/1}+C^{l(1)}_{i(m)}P^{(m)\;(1)}_
{(1)1(k)},}
$\medskip

$
\displaystyle{P^{l\;(1)}_{ij(k)}={\partial L^l_{ij}\over\partial
y^k}-C^{l(1)}_{i(k)\vert j}+C^{l(1)}_{i(m)}P^{(m)\;(1)}_
{(1)j(k)},}
$\medskip

$
\displaystyle{S^{l(1)(1)}_{i(j)(k)}={\partial C^{l(1)}_{i(j)}
\over\partial y^k}-{\partial C^{l(1)}_{i(k)}\over\partial y^j}
+C^{m(1)}_{i(j)}C^{l(1)}_{m(k)}-C^{m(1)}_{i(k)}C^{l(1)}_{m(j)};}
$\medskip

ii) for $p=\dim T\geq 2$, we have\medskip

$
\displaystyle{H^\alpha_{\eta\beta\gamma}={\partial H^\alpha_{\eta\beta}\over
\partial t^\gamma}-{\partial H^\alpha_{\eta\gamma}\over\partial t^\beta}+
H^\mu_{\eta\beta}H^\alpha_{\mu\gamma}-H^\mu_{\eta\gamma}H^\alpha_{\mu\beta},}
$\medskip

$
\displaystyle{R^l_{i\beta\gamma}={\delta G^l_{i\beta}\over\delta t^\gamma}-
{\delta G^l_{i\gamma}\over\delta t^\beta}+G^m_{i\beta}G^l_{m\gamma}-G^m_
{i\gamma}G^l_{m\beta},}
$\medskip

$
\displaystyle{R^l_{i\beta k}={\delta G^l_{i\beta}\over\delta x^k}-
{\delta\Gamma^l_{ik}\over\delta t^\beta}+G^m_{i\beta}\Gamma^l_{mk}-\Gamma^m_
{ik}G^l_{m\beta},}
$\medskip

$
\displaystyle{R^l_{ijk}=r^l_{ijk}={\partial\Gamma^l_{ij}\over\partial x^k}-
{\partial\Gamma^l_{ik}\over\partial x^j}+\Gamma^m_{ij}\Gamma^l_{mk}-\Gamma^
m_{ik}\Gamma^l_{mj}}
$.
\end{th}
\addtocounter{rem}{1}
{\bf Remark \therem} In the case of an autonomous metrical multi-time Lagrange
space of electrodynamics
(i. e. , $g_{ij}(t^\gamma,x^k,x^k_\gamma)=g_{ij}(x^k)$), all curvature d-tensors
of the Cartan canonical connection vanish,  except $H^\alpha_{\eta\beta\gamma}$
and\linebreak
$R^l_{ijk}=r^l_{ijk}$, that is, the curvature tensors  of the semi-Riemannian
metrics $h_{\alpha\beta}$ and  $g_{ij}$.

\section{Electromagnetic field. Maxwell equations}

\hspace{5mm} Let $ML^n_p=(J^1(T,M),L)$ be a metrical multi-time Lagrange
space and $\Gamma=(M^{(i)}_{(\alpha)\beta},N^{(i)}_{(\alpha)j})$ its canonical
nonlinear connection. Let us consider $C\Gamma=(H^\gamma_{\alpha\beta},G^k_{i\gamma},
L^k_{ij},C^{k(\gamma)}_{i(j)})$ the Cartan canonical connection of $ML^n_p$.

Using the {\it canonical Liouville d-tensor}  {\bf C}$\displaystyle{=x^i_
\alpha{\partial\over\partial x^i_\alpha}}$ and the fundamental vertical metrical
d-tensor $G^{(\alpha)(\beta)}_{(i)(k)}$ of the metrical multi-time Lagrange
space $ML^n_p$, we construct the {\it metrical deflection d-tensors}
\begin{equation}
\begin{array}{l}\medskip
\bar D^{(\alpha)}_{(i)\beta}=G^{(\alpha)(\gamma)}_{(i)(k)}\bar D^{(k)}_{
(\gamma)\beta}=x^{(\alpha)}_{(i)/\beta},\\\medskip
D^{(\alpha)}_{(i)j}=G^{(\alpha)(\gamma)}_{(i)(k)}D^{(k)}_{(\gamma)j}=
x^{(\alpha)}_{(i)\vert j},\\\medskip
d^{(\alpha)(\beta)}_{(i)(j)}=G^{(\alpha)(\gamma)}_{(i)(k)}d^{(k)(\beta)}_
{(\gamma)(j)}=x^{(\alpha)}_{(i)}\vert^{(\beta)}_{(j)},
\end{array}
\end{equation}
where $x^{(\alpha)}_{(i)}=G^{(\alpha)(\gamma)}_{(i)(k)}x^k_\gamma$ and
$"_{/\beta}"$, $"_{\vert j}"$ and $"\vert^{(\beta)}_{(j)}"$ are the local
covariant derivatives induced by $C\Gamma$.

Taking into account the expressions  of  the local covariant
derivatives of $C\Gamma$ (see the papers \cite{10}, \cite{13}), by  a direct  calculation, we obtain
\begin{prop}
The metrical deflection d-tensors of a metrical multi-time Lagrange space
$ML^n_p$ have the expressions:

i) for $p=1$,
\begin{equation}
\begin{array}{l}\medskip
\displaystyle{\bar D^{(1)}_{(i)1}={h^{11}\over 2}{\delta g_{im}\over\delta t}
y^m,}\\\medskip
D^{(1)}_{(i)j}=h^{11}g_{ik}\left[-N^{(k)}_{(1)j}+L^k_{jm}y^m\right],\\
d^{(1)(1)}_{(i)(j)}=h^{11}\left[g_{ij}+g_{ik}C^{k(1)}_{m(j)}y^m\right];
\end{array}
\end{equation}

ii) for $p\geq 2$,
\begin{equation}
\begin{array}{l}\medskip
\displaystyle{\bar D^{(\alpha)}_{(i)\beta}={h^{\alpha\gamma}\over 2}{\partial g_{km}\over
\partial t^\beta}x^m_\gamma,}\\\medskip
\displaystyle{D^{(\alpha)}_{(i)j}=-{h^{\alpha\gamma}\over 2}{\partial g_{ij}\over\partial
t^\gamma}-{1\over 4}U^{(\alpha)}_{(i)j},}\\
d^{(i)(\beta)}_{(\alpha)(j)}=h^{\alpha\beta}g_{ij}.
\end{array}
\end{equation}
\end{prop}
\addtocounter{rem}{1}
{\bf Remark \therem} In the particular case of an autonomous metrical multi-time
Lagrange space of electrodynamics
(i. e., $g_{ij}=g_{ij}(x^k)$), we have
$$
\bar D^{(\alpha)}_{(i)\beta}=0,\quad
D^{(\alpha)}_{(i)j}=-{1\over 4}U^{(\alpha)}_{(i)j},\quad
d^{(\alpha)(\beta)}_{(i)(j)}=h^{\alpha\beta}g_{ij}.
$$

In order to construct the metrical multi-time Lagrangian theory of electromagnetism,
we introduce the following\medskip\\
\addtocounter{defin}{1}
{\bf Definition \thedefin} The distinguished 2-form on $J^1(T,M)$,
\begin{equation}
F=F^{(\alpha)}_{(i)j}\delta x^i_\alpha\wedge dx^i+f^{(\alpha)(\beta)}_{(i)(j)}
\delta x^i_\alpha\wedge\delta x^j_\beta,
\end{equation}
where
$F^{(\alpha)}_{(i)j}=\displaystyle{{1\over 2}\left[D^{(\alpha)}
_{(i)j}-D^{(\alpha)}_{(j)i}\right]}$ and
$f^{(\alpha)(\beta)}_{(i)(j)}=\displaystyle{{1\over 2}\left[d^{(\alpha)(\beta)}
_{(i)(j)}-d^{(\alpha)(\beta)}_{(j)(i)}\right]}$, is called the \medskip\linebreak
{\it electromagnetic d-form} of the metrical multi-time Lagrange space
$ML^n_p$.\medskip\\
\addtocounter{rem}{1}
{\bf Remark \therem} The naturalness of the previous definition comes considering
the particular case of a relativistic rheonomic Lagrange space (i. e., $(T,h)=
(R,\delta)$). In this case, we recover the electromagnetic d-tensor of the Miron-Anastasiei
electromagnetism \cite{7}.\medskip

By simple computations, we find

\begin{prop}
The components $F^{(\alpha)}_{(i)j}$ and $f^{(\alpha)(\beta)}_{(i)(j)}$ of the
electromagnetic  d-form $F$ of a  metrical  multi-time Lagrange space are described
by the formulas:

i) in the case $p=1$,
$$
F^{(1)}_{(i)j}={h^{11}\over 2}\left[g_{jm}N^{(m)}_{(1)i}-g_{im}N^{(m)}_{(1)j}
+(g_{ik}L^k_{jm}-g_{jk}L^k_{im})y^m\right],\quad
f^{(1)(1)}_{(i)(j)}=0;
$$

ii) in the case $p\geq 2$,
$$
F^{(\alpha)}_{(i)j}={1\over 8}\left[U^{(\alpha)}_{(j)i}-U^{(\alpha)}_{(i)j}
\right],\quad f^{(\alpha)(\beta)}_{(i)(j)}=0.
$$
\end{prop}
\addtocounter{rem}{1}
{\bf Remark \therem} We emphasize that, in the particular case of an autonomous
metrical multi-time Lagrange space (i. e. $g_{ij}=g_{ij}(x^k)$), the electromagnetic
components get the expressions
$$
F^{(\alpha)}_{(i)j}={1\over 8}\left[U^{(\alpha)}_{(j)i}-U^{(\alpha)}_{(i)j}
\right],\quad f^{(\alpha)(\beta)}_{(i)(j)}=0
$$

The main result of the electromagnetic metrical  multi-time Lagrangian theory
is the following
\begin{th}
The electromagnetic components $F^{(\alpha)}_{(i)j}$ of a metrical multi-time
Lagrange space $ML^n_p=(J^1(T,M),L)$ are governed by the Maxwell equations:
\medskip

i) for $p=1$,
$$
\left\{\begin{array}{l}\medskip
\displaystyle{F^{(1)}_{(i)k/1}={1\over 2}{\cal A}_{\{i,k\}}\left\{\bar D^{(1)}_{(i)1\vert k}
+D^{(1)}_{(i)m}T^m_{1k}+d^{(1)(1)}_{(i)(m)}R^{(m)}_{(1)1k}-\left[T^p_{1i\vert k}
+C^{p(1)}_{k(m)}R^{(m)}_{(1)1i}\right] y_{(p)}\right\}}\\\medskip
\displaystyle{\sum_{\{i,j,k\}}F^{(1)}_{(i)j\vert k}=-{1\over 2}\sum_{\{i,j,k\}}C^{(1)(1)(1)}_
{(i)(l)(m)}R^{(m)}_{(1)jk}y^l}\\
\displaystyle{\sum_{\{i,j,k\}}F^{(1)}_{(i)j}\vert^{(1)}_{(k)}=0,}
\end{array}\right.
$$

ii) for $p\geq 2$,
$$
\left\{\begin{array}{l}\medskip
\displaystyle{F^{(\alpha)}_{(i)k/\beta}={1\over 2}{\cal A}_
{\{i,k\}}\left\{\bar D^{(\alpha)}_{(i)\beta\vert k}
+D^{(\alpha)}_{(i)m}T^m_{\beta k}+d^{(\alpha)(\mu)}_
{(i)(m)}R^{(m)}_{(\mu)\beta k}-\left[T^p_{\beta i\vert k}+C^{p(\mu)}_
{k(m)}R^{(m)}_{(\mu)\beta i}\right]
x^{(\alpha)}_{(p)}\right\}}\\\medskip
\displaystyle{\sum_{\{i,j,k\}}F^{(\alpha)}_{(i)j\vert k}=0}\\
\displaystyle{\sum_{\{i,j,k\}}F^{(\alpha)}_{(i)j}\vert^{(\gamma)}_{(k)}=0,}
\end{array}\right.
$$
where $y_{(p)}=G^{(1)(1)}_{(p)(q)}y^q$,  $\displaystyle{C^{(1)(1)(1)}_{(i)(l)(m)}=G^{(1)(1)}_{(l)(q)}C^{q(1)}_{i(m)}
={h^{11}\over 2}{\partial^3L\over\partial y^i\partial y^l\partial y^m}}$,
$x^{(\alpha)}_{(p)}=G^{(\alpha)(\beta)}_{(p)(q)}x^q_\beta$.
\end{th}
{\bf Proof.} Firstly, we point out that the Ricci identities \cite{13} applied
to the spatial metrical d-tensor $g_{ij}$ imply that the following curvature
d-tensor identities
$$
R_{mi\beta k}+R_{im\beta k}=0,\quad
R_{mijk}+R_{imjk}=0,\quad
P^{\;\;\;\;\;(\gamma)}_{mij(k)}+P^{\;\;\;\;\;(\gamma)}_{imj(k)}=0,
$$
where $R_{mi\beta k}=g_{ip}R^p_{m\beta k}$,
$R_{mijk}=g_{ip}R^p_{mjk}$ and
$P^{\;\;\;\;\;(\gamma)}_{mij(k)}=g_{ip}P^{p\;\;\;(\gamma)}_{mj(k)}$, are true.

Now, let us consider the following general deflection d-tensor identities
\cite{13}\medskip

$d_1)\;\;\bar D^{(p)}_{(\nu)\beta\vert k}-D^{(p)}_{(\nu)k/\beta}=x^m_\nu
R^p_{m\beta k}-D^{(p)}_{(\nu)m}T^m_{\beta k}-d^{(p)(\mu)}_{(\nu)(m)}R^{(m)}_
{(\mu)\beta k},$\medskip

$d_2)\;\;D^{(p)}_{(\nu)j\vert k}-D^{(p)}_{(\nu)k\vert j}=x^m_\nu
R^p_{mjk}-d^{(p)(\mu)}_{(\nu)(m)}R^{(m)}_{(\mu)jk},$\medskip

$d_3)\;\;D^{(p)}_{(\nu)j}\vert^{(\gamma)}_{(k)}-d^{(p)(\gamma)}_{(\nu)(k)
\vert j}=x^m_\nu P^{p\;\;(\gamma)}_{mj(k)}-D^{(p)}_{(\nu)m}C^{m(\gamma)}_
{j(k)}-d^{(p)(\mu)}_{(\nu)(m)}P^{(m)\;(\gamma)}_{(\mu)j(k)}$,\\
where
$\bar D^{(i)}_{(\alpha)\beta}=x^i_{\alpha/\beta}$,
$D^{(i)}_{(\alpha)j}=x^i_{\alpha\vert j}$,
$d^{(i)(\beta)}_{(\alpha)(j)}=x^i_\alpha\vert^{(\beta)}_{(j)}$.
Contracting the deflection d-tensor identities by $G^{(\alpha)(\nu)}_{(i)(p)}$
and using the above curvature d-tensor equalities, we obtain the metrical deflection
d-tensors  identities:\medskip

$d^\prime_1)\;\;\bar D^{(\alpha)}_{(i)\beta\vert k}-D^{(\alpha)}_{(i)k/\beta}=
-x^{(\alpha)}_{(m)}R^m_{i\beta k}-D^{(\alpha)}_{(i)m}T^m_{\beta k}-d^{(\alpha)(\mu)}
_{(i)(m)}R^{(m)}_{(\mu)\beta k},$\medskip

$d^\prime_2)\;\;D^{(\alpha)}_{(i)j\vert k}-D^{(\alpha)}_{(i)k\vert j}=
-x^{(\alpha)}_{(m)}R^m_{ijk}-d^{(\alpha)(\mu)}_{(i)(m)}R^{(m)}_{(\mu)jk},
$\medskip

$d^\prime_3)\;\;D^{(\alpha)}_{(i)j}\vert^{(\gamma)}_{(k)}-d^{(\alpha)(\gamma)}
_{(i)(k)\vert j}=-x^{(\alpha)}_{(m)}P^{m\;\;(\gamma)}_{ij(k)}-D^{(\alpha)}
_{(i)m}C^{m(\gamma)}_{j(k)}-d^{(\alpha)(\mu)}_{(i)(m)}P^{(m)\;(\gamma)}_
{(\mu)j(k)}$.\medskip

At the same time, we recall that the following Bianchi identities \cite{10}\medskip

$b_1)\;\;{\cal A}_{\{j,k\}}\left\{R^l_{j\alpha k}+T^l_{\alpha j\vert k}+
C^{l(\mu)}_{k(m)}R^{(m)}_{(\mu)\alpha j}\right\}=0$,\medskip

$b_2)\;\;\sum_{\{i,j,k\}}\left\{R^l_{ijk}-
C^{l(\mu)}_{k(m)}R^{(m)}_{(\mu)ij}\right\}=0$,\medskip

$b_3)\;\;{\cal A}_{\{j,k\}}\left\{P^{l\;\;(\varepsilon)}_{jk(p)}+C^{l(\varepsilon)}
_{j(p)\vert k}+C^{l(\mu)}_{k(m)}P^{(m)\;\;(\varepsilon)}_{(\mu)j(p)}\right\}=0$,
\medskip\\
where ${\cal A}_{\{j,k\}}$ means alternate sum and $\sum_{\{i,j,k\}}$ means
cyclic sum, hold good.

In order to obtain the first Maxwell identity, we permute $i$ and $k$ in
$d^\prime_1$ and we subtract the new identity from the initial one. Finally,
using the Bianchi identity $b_1$, we obtain what we were looking for.

Doing a cyclic sum by the indices $\{i,j,k\}$ in $d^\prime_2$ and using the
Bianchi identity $b_2$, it follows the second Maxwell equation.

Applying a Christoffel process to the indices $\{i,j,k\}$ in $d^\prime_3$
and combining with the Bianchi identity $b_3$ and the relation
$P^{(m)\;\;(\varepsilon)}_{(\mu)j(p)}=P^{(m)\;\;(\varepsilon)}_{(\mu)p(j)}$,
we get a new identity. The cyclic sum by the indices
$\{i,j,k\}$ applied to this last identity implies the third Maxwell equation.
\rule{5pt}{5pt}\medskip\\
\addtocounter{rem}{1}
{\bf Remark \therem} In the case of an autonomous metrical multi-time Lagrange
space of electrodynamics (i. e., $g_{ij}=g_{ij}(x^k)$), the Maxwell equations
take a more simple form, namely,
$$
F^{(\alpha)}_{(i)k/\beta}={1\over 2}{\cal A}_{\{i,k\}}
h^{\alpha\mu}g_{im}R^{(m)}_{(\mu)\beta k},\quad
\sum_{\{i,j,k\}}F^{(\alpha)}_{(i)j\vert k}=0,\quad
\sum_{\{i,j,k\}}F^{(\alpha)}_{(i)j}\vert^{(\gamma)}_{(k)}=0.
$$

\section{Gravitational field. Einstein equations}

\setcounter{equation}{0}
\hspace{5mm} Let $h=(h_{\alpha\beta})$ be a fixed semi-Riemannian metric on
the temporal manifold $T$ and $\Gamma=(M^{(i)}_{(\alpha)\beta},
N^{(i)}_{(\alpha)j})$ a fixed nonlinear connection on the 1-jet space
$J^1(T,M)$. In order to develope the metrical  multi-time Lagrange theory
of gravitational field, we introduce the following\medskip\\
\addtocounter{defin}{1}
{\bf Definition \thedefin} From physical point of view, an adapted metrical d-tensor
$G$ on $E=J^1(T,M)$, expressed locally by
$$
G=h_{\alpha\beta}dt^\alpha\otimes dt^\beta+g_{ij}dx^i\otimes dx^j+h^{\alpha
\beta}g_{ij}\delta x^i_\alpha\otimes\delta x^j_\beta,
$$
where $g_{ij}=g_{ij}(t^\gamma,x^k,x^k_\gamma)$ is a d-tensor field on $E$,
symmetric, of rank $n=\dim M$ and having a constant signature on $E$, is
called a {\it gravitational $h$-potential}.\medskip\\
\addtocounter{rem}{1}
{\bf Remark \therem} The naturalness of this definition comes from the particular
case $(T,h)=(R,\delta)$. In this case, we recover the {\it gravitational
potentials} $g_{ij}(x,y)$ from Miron-Anastasiei theory of gravitational field
\cite{7}.\medskip

Now, taking $ML^n_p=(J^1(T,M),L)$ a metrical multi-time Lagrange space, via
its fundamental vertical metrical d-tensor
\begin{equation}
G^{(\alpha)(\beta)}_{(i)(j)}={1\over 2}{\partial^2L\over\partial x^i_\alpha
\partial x^j_\beta}=\left\{\begin{array}{ll}\medskip
h^{11}(t)g_{ij}(t,x^k,y^k),&p=\dim T=1\\
h^{\alpha\beta}(t^\gamma)g_{ij}(t^\gamma,x^k),&p=\dim T\geq 2,
\end{array}\right.
\end{equation}
and its canonical nonlinear connection $\Gamma=(M^{(i)}_{(\alpha)\beta},
N^{(i)}_{(\alpha)j})$, one induces a natural gravitational $h$-potential,
setting
$$
G=h_{\alpha\beta}dt^\alpha\otimes dt^\beta+g_{ij}dx^i\otimes dx^j+h^{\alpha
\beta}g_{ij}\delta x^i_\alpha\otimes\delta x^j_\beta.
$$
Let us consider $C\Gamma=(H^\gamma_{\alpha\beta},G^k_{j\gamma},L^i_{jk},
C^{i(\gamma)}_{j(k)})$ the Cartan canonical connection of $ML^n_p$.

We postulate that the Einstein which govern the gravitational $h$-potential
$G$ of the metrical multi-time Lagrange space $ML^n_p$ are the Einstein equations attached
to the Cartan canonical connection $C\Gamma$ of $ML^n_p$ and the adapted metric $G$ on
$E$, that is,
\begin{equation}
Ric(C\Gamma)-{Sc(C\Gamma)\over 2}G={\cal K}{\cal T},
\end{equation}
where $Ric(C\Gamma)$ represents the Ricci d-tensor of the Cartan connection,
$Sc(C\Gamma)$ is its scalar curvature, ${\cal K}$ is the Einstein constant and ${\cal T}$
is an intrinsec tensor of matter which is called  the {\it stress-energy}
d-tensor.

In the adapted basis $(X_A)=\displaystyle{\left({\delta\over\delta t^\alpha},
{\delta\over\delta x^i},{\partial\over\partial x^i_\alpha}\right)}$ of  the
nonlinear  connection $\Gamma$ of $ML^n_p$, the curvature d-tensor {\bf R} of
the Cartan connection is expressed locally by {\bf R}$(X_C,X_B)X_A=R^D_{ABC}X_D$.
It follows that we have $R_{AB}=Ric(X_A,X_B)=R^D_{ABD}$ and $Sc(C\Gamma)
=G^{AB}R_{AB}$, where
\begin{equation}
G^{AB}=\left\{\begin{array}{ll}\medskip
h_{\alpha\beta},&\mbox{for}\;\;A=\alpha,\;B=\beta\\\medskip
g^{ij},&\mbox{for}\;\;A=i,\;B=j\\\medskip
h_{\alpha\beta}g^{ij},&\mbox{for}\;\;A={(i)\atop(\alpha)},\;B={(j)\atop(\beta)}\\
0,&\mbox{otherwise}.
\end{array}\right.
\end{equation}

Taking into account, on the one hand, the form of the fundamental vertical
metrical  d-tensor  $G^{(\alpha)(\beta)}_{(i)(j)}$   of the metrical multi-time
Lagrange space $ML^n_p$, and,  on the other  hand, the expressions of local
curvature d-tensors attached to the  Cartan canonical  connection $C\Gamma$,
by a direct calculation, we deduce
\begin{th}
The Ricci d-tensor $Ric(C\Gamma)$ of the  Cartan canonical connection $C\Gamma$
of a metrical multi-time Lagrange space, is determined by the following components:
\medskip

i) for $p=\dim T=1$,
$$
\begin{array}{l}\medskip
R_{11}\stackrel{not}{=}H_{11}=0,\quad R_{i1}=R^m_{i1m},\quad R_{ij}=R^m_{ijm},\quad
R^{\;(1)}_{i(j)}\stackrel{not}{=}P^{\;(1)}_{i(j)}=-P^{m\;\;(1)}_{im(j)},\\
R^{(1)}_{(i)j}\stackrel{not}{=}P^{(1)}_{(i)j}=P^{m\;(1)}_{ij(m)},\quad
R^{(1)}_{(i)1}\stackrel{not}{=}P^{(1)}_{(i)1}=P^{m\;(1)}_{i1(m)},\quad
R^{(1)(1)}_{(i)(j)}\stackrel{not}{=}S^{(1)(1)}_{(i)(j)}=S^{m(1)(1)}_{i(j)(m)};
\end{array}
$$

ii) for $p=\dim T\geq 2$,
$$
\begin{array}{l}\medskip
R_{(\alpha)(\beta)}\stackrel{not}{=}H_{\alpha\beta}=H^\mu_{\alpha\beta\mu},\quad
R_{i\alpha}=R^m_{i\alpha m},\quad R_{ij}=R^m_{ijm},\quad
R^{\;(\alpha)}_{i(j)}\stackrel{not}{=}P^{\;(\alpha)}_{i(j)}=0,\\
R^{(\alpha)}_{(i)j}\stackrel{not}{=}P^{(\alpha)}_{(i)j}=0,\quad
R^{(\alpha)}_{(i)\beta}\stackrel{not}{=}P^{(\alpha)}_{(i)\beta}=0,\quad
R^{(\alpha)(\beta)}_{(i)(j)}\stackrel{not}{=}S^{(\alpha)(\beta)}_{(i)(j)}=0.
\end{array}
$$
\end{th}

Denoting $H=h^{\alpha\beta}H_{\alpha\beta},\;R=g^{ij}R_{ij}$ and
$S=h_{\alpha\beta}g^{ij}S^{(\alpha)\beta)}_{(i)(j)}$, it follows
\begin{cor}
The scalar curvature of $Sc(C\Gamma)$ of the Cartan canonical connection $C\Gamma$
of   a  metrical multi-timeLagrange space, is given  by the formulas
\medskip

i) for $p=\dim T=1,\qquad Sc(C\Gamma)=R+S;$\medskip

ii) for $p=\dim T\geq 2,\qquad Sc(C\Gamma)=H+R.$
\end{cor}
\addtocounter{rem}{1}
{\bf Remark \therem} In the particular case of an autonomous metrical multi-time
Lagrange space of electrodynamics (i. e., $g_{ij}=g_{ij}(x^k)$), all Ricci
d-tensor components vanish, except $H_{\alpha\beta}$ and $R_{ij}=r_{ij}$,
where $H_{\alpha\beta}$ (resp. $r_{ij}$) are the local Ricci tensors associated
to the semi-Riemannian metric $h_{\alpha\beta}$ (resp. $g_{ij}$). It follows
that the scalar curvature of a this space is $Sc(C\Gamma)=H+r$, where $H$ and
$r$ are the scalar curvatures of the semi-Riemannian metrics $h_{\alpha\beta}$
and $g_{ij}$.\medskip

The main result of the  metrical multi-time Lagrange theory of gravitational
field is  given by the  following
\begin{th}
The Einstein equations which govern the gravitational $h$-potential
$G$ induced by the Kronecker $h$-regular Lagrangian of a metrical multi-time
Lagrange space $ML^n_p$, take the form\medskip

i) for $p=\dim T=1$,
$$
\left\{\begin{array}{l}\medskip
\displaystyle{-{R+S\over 2}h_{11}={\cal K}{\cal T}_{11}}\\\medskip
\displaystyle{R_{ij}-{R+S\over 2}g_{ij}={\cal K}{\cal T}_{ij}}\\\medskip
\displaystyle{S^{(1)(1)}_{(i)(j)}-{R+S\over 2}h^{11}g_{ij}={\cal K}{\cal T}
^{(1)(1)}_{(i)(j)}},
\end{array}\right.\leqno{(E_1)}
$$
$$
\left\{\begin{array}{lll}\medskip
0={\cal T}_{1i},&R_{i1}={\cal K}{\cal T}_{i1},&
P^{(1)}_{(i)1}={\cal K}{\cal T}^{(1)}_{(i)1}\\
0={\cal T}^{\;(1)}_{1(i)},&P^{\;(1)}_{i(j)}={\cal K}{\cal T}^{\;(1)}_{i(j)},&
P^{(1)}_{(i)j}={\cal K}{\cal T}^{(1)}_{(i)j},
\end{array}\right.\leqno{(E_2)}
$$

i) for $p=\dim T\geq 2$,
$$
\left\{\begin{array}{l}\medskip
\displaystyle{H_{\alpha\beta}-{H+R\over 2}h_{\alpha\beta}={\cal K}{\cal T}_
{\alpha\beta}}\\\medskip
\displaystyle{R_{ij}-{H+R\over 2}g_{ij}={\cal K}{\cal T}_{ij}}\\\medskip
\displaystyle{-{H+R\over 2}h^{\alpha\beta}g_{ij}={\cal K}{\cal T}^{(\alpha)
(\beta)}_{(i)(j)}},
\end{array}\right.\leqno{(E_1)}
$$
$$
\left\{\begin{array}{lll}\medskip
0={\cal T}_{\alpha i},&R_{i\alpha}={\cal K}{\cal T}_{i\alpha},&
0={\cal T}^{(\alpha)}_{(i)\beta}\\
0={\cal T}^{\;(\beta)}_{\alpha(i)},&0={\cal T}^{\;(\alpha)}_{i(j)},&
0={\cal T}^{(\alpha)}_{(i)j},
\end{array}\right.\leqno{(E_2)}
$$
where ${\cal T}_{AB},\;A,B\in\{\alpha,i,{(\alpha)\atop(i)}\}$ are the adapted
local components of the stress-energy d-tensor ${\cal T}$.
\end{th}
\addtocounter{rem}{1}
{\bf Remarks \therem} i) Asumming that $p=\dim T>2$ and $n=\dim M>2$, the set
$(E_1)$ of the Einstein equations can be rewritten in the more natural form
$$
\left\{\begin{array}{l}\medskip
\displaystyle{H_{\alpha\beta}-{H\over 2}h_{\alpha\beta}={\cal K}\tilde{\cal T}_
{\alpha\beta}}\\\medskip
\displaystyle{R_{ij}-{R\over 2}g_{ij}={\cal K}\tilde{\cal T}_{ij}},
\end{array}\right.\leqno{(E^\prime_1)}
$$
where $\tilde{\cal T}_{AB},\;A,B\in\{\alpha,i\}$ are the adapted local components
of a new stress-energy d-tensor $\tilde{\cal T}$. This new form of the Einstein
equations will be treated detalied in the more general case of a {\it generalized
metrical multi-time Lagrange space} \cite{9}.

ii) In the particular case of an autonomous metrical multi-time Lagrange space
of electrodynamics (i. e., $g_{ij}=g_{ij}(x^k)$), the following
Einstein equations of gravitational field
$$
\left\{\begin{array}{l}\medskip
\displaystyle{H_{\alpha\beta}-{H+r\over 2}h_{\alpha\beta}={\cal K}{\cal T}_
{\alpha\beta}}\\\medskip
\displaystyle{r_{ij}-{H+r\over 2}g_{ij}={\cal K}{\cal T}_{ij}}\\\medskip
\displaystyle{-{H+r\over 2}h^{\alpha\beta}g_{ij}={\cal K}{\cal T}^{(\alpha)
(\beta)}_{(i)(j)}},
\end{array}\right.\leqno{(E_1)}
$$
$$
\left\{\begin{array}{lll}\medskip
0={\cal T}_{\alpha i},&0={\cal T}_{i\alpha},&
0={\cal T}^{(\alpha)}_{(i)\beta}\\
0={\cal T}^{\;(\beta)}_{\alpha(i)},&0={\cal T}^{\;(\alpha)}_{i(j)},&
0={\cal T}^{(\alpha)}_{(i)j},
\end{array}\right.\leqno{(E_2)}
$$
hold good. It is remarkable that the new form $(E_1^\prime)$ of the Einstein
equations of a metrical multi-time Lagrange space of electrodynamics reduces
to the classical one, namely,
$$
\left\{\begin{array}{l}\medskip
\displaystyle{H_{\alpha\beta}-{H\over 2}h_{\alpha\beta}={\cal K}\tilde{\cal T}_
{\alpha\beta}}\\\medskip
\displaystyle{r_{ij}-{r\over 2}g_{ij}={\cal K}\tilde{\cal T}_{ij}}.
\end{array}\right.
$$

iii) In order to have the compatibility of the Einstein equations, it is
necessary that the certain  adapted local components of the stress-energy
d-tensor vanish {\it "a priori"}.\medskip

From physical point of view, it is  well known that the stress-energy
d-tensor ${\cal T}$ must  verify the local {\it conservation laws} ${\cal T}^
B_{A\vert B}=0,\;\forall\;A\in\{\alpha,i,{(\alpha)\atop (i)}\}$,
where ${\cal T}^B_A=G^{BD}{\cal T}_{DA}$. Consequently, by a direct calculation,
we find the  following
\begin{th}
The conservation laws of the Einstein equations of a metrical multi-time
Lagrange space $ML^n_p$ are given by the  formulas\medskip

i) for $p=1$,
\begin{equation}
\left\{\begin{array}{l}\medskip
\displaystyle{\left[{R+S\over 2}\right]_{/1}=R^m_{1\vert m}-P^{(m)}_{(1)1}
\vert^{(1)}_{(m)}}\\\medskip
\displaystyle{\left[R^m_j-{R+S\over 2}\delta^m_j\right]_{\vert m}=
-P^{(m)}_{(1)j}\vert^{(1)}_{(m)}}\\\medskip
\displaystyle{\left[S^{(m)(1)}_{(1)(j)}-{R+S\over 2}\delta^m_j\right]\left\vert
^{(1)}_{(m)}\right.=-P^{m(1)}_{\;\;\;(j)\vert m}},
\end{array}\right.
\end{equation}\medskip
where
$
R^i_1=g^{im}R_{m1},\;\; P^{(i)}_{(1)1}=h_{11}g^{im}P^{(1)}_{(m)1},\;\;
R^i_j=g^{im}R_{mj},\;\;
P^{(i)}_{(1)j}=h_{11}g^{im}P^{(1)}_{(m)j},
$
$
P^{i(1)}_{\;\;(j)}=g^{im}P^{\;\;(1)}_{m(j)}\;\mbox{and}\;
S^{(i)(1)}_{(1)(j)}=h_{11}g^{im}S^{(1)(1)}_{(m)(j)};
$\medskip

i) for $p\geq 2$,
\begin{equation}
\left\{\begin{array}{l}\medskip
\displaystyle{\left[H^\mu_\beta-{H+R\over 2}\delta^\mu_\beta\right]_{/\mu}=
-R^m_{\beta\vert m}}\\
\displaystyle{\left[R^m_j-{H+R\over 2}\delta^m_j\right]_{\vert m}=0},
\end{array}\right.
\end{equation}
where
$
H^\mu_\beta=h^{\mu\gamma}H_{\gamma\mu},\;\;
R^i_j=g^{im}R_{mj}\;\mbox{and}\;R^i_\beta=g^{im}R_{m\beta}.
$
\end{th}
\addtocounter{rem}{1}
{\bf Remarks \therem} i) In the case $p>2,\;n>2$, taking into account the components
$\tilde{\cal T}_{\alpha\beta}$ and $\tilde{\cal T}_{ij}$ of the new stress-energy
d-tensor ${\cal T}$ from $(E^\prime_1)$, we point out that the conservation
laws modify in the following simple and natural  new form
$\tilde{\cal T}^\mu_{\beta/\mu}=0,\;\tilde{\cal T}^m_{j\vert m}=0$ (see \cite{9}).

ii) Considering an autonomous metrical multi-time Lagrange space
of electrodynamics (i. e., $g_{ij}=g_{ij}(x^k)$), the conservation laws of the
Einstein equations reduce to
$$
\left\{\begin{array}{l}\medskip
\displaystyle{\left[H^\mu_\beta-{H+r\over 2}\delta^\mu_\beta\right]_{/\mu}=0}
\\
\displaystyle{\left[r^m_j-{H+r\over 2}\delta^m_j\right]_{\vert m}=0}.
\end{array}\right.
$$

\section{Conclusion}

\hspace{5mm} Note that all entities with geometrical or physical meaning from this paper
was directly arised from the fundamental vertical metrical d-tensor $G^{(\alpha)(\beta)}_
{(i)(j)}$ of $ML^n_p$.  This fact points out the {\it metrical character}
and the naturalness of the metrical multi-time Lagrange theory of physical
fields that we constructed. At the same  time, the form of the invariance
gauge group \ref{gg} of the fibre bundle of configurations, $J^1(T,M)\to T\times M$,
allows us to appreciate the metrical multi-time Lagrangian field theory like
a {\it "parametrized"} theory. In conclusion, the metrical multi-time Lagrangian
theory of physical fields is, via the Marsden's classification of field
theories \cite{4}, a {\it "metrical-parametrized"} one.\medskip\\
{\bf\underline{Open problem}.} The development of an analogous metrical
multi-time Lagrangian geometry
of physical fields on the jet space of order two $J^2(T,M)$  is in our attention.
\medskip\\
{\bf Acknowledgements.} A version of this paper was presented at
Workshop on Diff. Geom. , Global Analysis, Lie Algebras, Aristotle University of Thessaloniki,
Greece, Aug. 28-Sept. 2, 2000. The author thanks to all participants at this
conference, especially to Prof. Dr. C. Udriste for its useful sugestions.

\begin{center}
University POLITEHNICA of Bucharest\\
Department of Mathematics I\\
Splaiul Independentei 313\\
77206 Bucharest, Romania\\
e-mail: mircea@mathem.pub.ro\\
\end{center}
\end{document}